\documentclass[12pt]{article}
\usepackage{amsmath,amssymb,amscd,latexsym,euscript,mathrsfs}
\usepackage[matrix,arrow,curve]{xy}

\newcommand{\Ob}{ {\rm Ob\,}}
\newcommand{\Mor}{ {\rm Mor\,}}
\newcommand{\dom}{ {\rm dom\,}}
\newcommand{\cod}{ {\rm cod\,}}
\newcommand{\coLim}{\underrightarrow{\lim}}
\newcommand{\Lim}{\underleftarrow{\lim}}
\newcommand{\gld}{{\rm gl.dim\,}}

\newcommand{\Set}{{\rm Set}}
\newcommand{\Ab}{{\rm Ab}}

\newcommand{\Imm}{{\rm Im\,}}
\newcommand{\Ker}{{\rm Ker\,}}

\newcommand{\Ext}{{\rm Ext\,}}

\newcommand{\NN}{{\,\mathbb N}}

\newcommand{\mC}{{\mathscr C}}

\newcommand{\mN}{{\,\cal N}}

\newcommand{\ZZ}{{\,\mathbb Z}}
\newcommand{\II}{{\,\mathbb I}}
\newcommand{\mA}{{\mathcal A}}

\newcommand{\mS}{{\mathcal S}}

\newtheorem{theorem}{\bf Theorem}[section]

\newtheorem{proposition}[theorem]{\bf Proposition}
\newtheorem{corollary}[theorem]{\bf Corollary}
\newtheorem{definition}{\sc Definition}[section]
\newtheorem{example}[definition]{\sc Example}

\newtheorem{conjecture}{\sc Conjecture}
\newtheorem{problem}{\sc Problem}

\def\leq{\leqslant}
\def\geq{\geqslant}

\begin{document}

\begin{center}
 {\large THE CUBICAL HOMOLOGY OF TRACE MONOIDS
}
\end{center}
\centerline{A. A. Husainov, husainov51@yandex.ru}

\begin{abstract}
This article contains an overview of the results of the author
in a field of algebraic topology used in computer science.
The relationship between the cubical homology groups of generalized tori
and  homology groups of partial trace monoid actions is described.
Algorithms for computing the homology groups of asynchronous systems, Petri nets,
and Mazurkiewicz trace languages are shown.
\end{abstract}


2000 Mathematics Subject Classification 18G10, 18G35, 55U10, 68Q10, 68Q85

\section*{Introduction}

Trace monoids have found many applications in computer science \cite{maz1987}, \cite{die1997}.
M. Bednarczyk \cite{bed1988} studied and applied the category of asynchronous systems.
The author has proved that any asynchronous system can be regarded as a partial
trace monoid with action on a set.
It is allowed to build homology theory for the category of asynchronous systems
and Petri nets \cite{X20042}. 
In this paper we also introduce a homology for Mazurkiewicz trace languages.
It should be noted that the homology theory was introduced 
and studied for higher dimensional automata in the works \cite{gou1995} and \cite{gau2000}.
E. Haucourt \cite{hau2009} applied the Baues-Wirsching homology.

We study of a relationship between the cubical homology
of generalized tori and homology of a trace monoid action on a set.
We build the algorithms for the computing the homology
groups of asynchronous systems, elementary Petri nets, and Mazurkiewicz trace languages.
It allows us to solve the problem posed in \cite[Open problem 1]{X20042}
of the constructing an algorithm for computing homology groups
of the elementary Petri nets.


\section{Trace monoids and their partial actions}

 This section is devoted to the basic definitions and
the problems that have appeared.

\subsection{Notations}

We first describe our notations.
Let $\Set$ be a category of sets and maps and let $\Ab$ be a category of Abelian groups
and homomorphisms. Denote by $\ZZ$ the set of the additive group of integers.
Let $\NN$ be the set of nonnegative integers or the free monoid $\{1, a, a^2, \cdots \}$
 generated by one element.
Given a category $\mA$ denote by $\mA^{op}$ the opposite category.
Denote by $\Ob\mA$ the class of all objects and $\Mor\mA$ the class of all morphisms in
the category $\mA$.
 Given objects $a,b\in \Ob\mA$ denote by
$\mA(a,b)$ the set of all morphisms $a\rightarrow b$.
For any small category $\mC$, functors $F: \mC\to \mA$
will be called {\em diagrams of objects in $\mA$ on $\mC$}.
In this case, along with the notation $F: \mC\to \mA$ we use the notation $\{F(c)\}_{c\in \mC}$.
The category $\mA^{\mC}$ of functors $\mC\to \mA$
is called to be a {\em diagram category}.

Let $\Delta\ZZ: \mC\to \Ab$ be the diagram having the value
$\ZZ$ at each $c\in Ob\mC$ and the value
 $1_\ZZ$ at each $\alpha\in Mor\mC$.

Given a family of Abelian groups $\{A_j\}_{j\in J}$ the direct sum is denoted
by $\bigoplus\limits_{j\in J}A_j$. Elements of summands are denoted as pairs
 $(j,g)$ with $j\in J$ and $g\in A_j$.
If $A_j= A$ for all $j\in J$, then this direct is denoted
 $\bigoplus\limits_{j\in J}A$. If instead of a set  $J$
indicated a cardinal number $p=|J|$, then the direct coproduct is denoted by $A^{(p)}$.


\subsection{Trace monoids}
Let $E$ be a set and let $I\subseteq E\times E$ be an arbitrary subset.
The set $I\subseteq E\times E$ is an {\em independence relation} on $E$ if the following
conditions are satisfied:

\begin{itemize}
\item $(\forall a\in E) (a,a)\notin I$~,
\item $(\forall a\in E)(\forall b\in E)~ (a,b)\in I \Rightarrow (b,a)\in I$.
\end{itemize}

Let $E^*$ be the free monoid generated by a set $E$.
It consists of the words in alphabet $E$. The binary operation is defined as the concatenation of words
$(a_1\cdots a_m, b_1\cdots b_n)\mapsto a_1\cdots a_m b_1\cdots b_n$.
The empty word is denoted by $1$.

\begin{definition}
Let $I$ be  an independence relation on a set $E$.
 A trace monoid (or free partially commutative monoid) $M(E,I)$ is the factor monoid $E^*/(\equiv)$
by a least equivalence relation for which $u a b v \equiv u b a v$, for all $(a,b)\in I$,
$u\in E^*$, $v\in E^*$.
Elements $a,b\in E$ for which $(a,b)\in I$ are called commuting generators.
\end{definition}
This definition is more general than given in \cite{die1997}
since we do not demand that $E$ is finite.

For example, if $E =\{a,b\}$, $I =\{(a,b),(b,a)\}$, then $M(E,I) \cong \NN^2$
is the free commutative monoid generated by two elements.

If $I=\emptyset$, then $M(E,I)=E^*$.

Any element $w=a_1\cdots a_n \in M(E,I)$ of a trace monoid can be interpreted as finite
sequence of instructions $a_1, a_2, \cdots, a_n$ in a program. The relation $I$ consists
of pairs $(a,b)$ instructions which can be executed concurrenrly.

\subsection{State space}

A {\em partial map} $f: E \rightharpoonup E'$ between sets $E$ and $E'$
is a relation $f \subseteq E\times E'$ for which $(e, e'_1)\in f ~\& ~(e, e'_2)\in f$
implies $e'_1=e'_2$.
Let $PSet$ be the category of sets and partial maps between them.
Any trace monoid $M(E,I)$ can be considered
as a category with the unique object denoted by $o(M(E,I))$.

A {\em partial trace monoid action} of $M(E,I)$ on a set $S$ is a functor
${\bf S}: M(E,I)^{op}\to PSet$ such that its value at $o(M(E,I))$
equals $S$. We denote ${\bf S}(w)(s)$ by $s\cdot w$.
A {\em state space} $(M(E,I),S)$ consists of a trace monoid
$M(E,I)$ with a partial action on a set $S$.
A state space $(M(E,I),S)$ is determined by partial maps $(-)\cdot a: S \rightharpoonup S$
corresponding to $a\in E$. Hence, it can be given by a directed graph with
vertexes $s\in S$ and labeled edges $s \stackrel{a}\to s\cdot e$.

For example, if $E=\{a,b\}$ and $I=\{(a,b), (b,a)\}$, then the directed graph
with labeled edges
$$
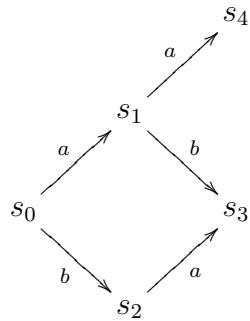
\xymatrix{
& & s_4\\
& s_1 \ar[ru]^{a}\ar[rd]^{b}\\
s_0 \ar[ru]^a \ar[rd]_b && s_3\\
& s_2 \ar[ru]_a
}
$$
determines the action for which $s_0\cdot a=s_1$, $s_0\cdot b=s_2$,
$s_1\cdot a=s_4$, $s_1\cdot b=s_3$,
$s_2\cdot a=s_3$. But $s_2\cdot b$, $s_3\cdot a$, $s_3\cdot b$,
$s_4\cdot a$, and $s_4\cdot b$ are not defined.

\subsection{Augmented state space}

In order to make the action $(M(E,I),S)$ to be total, we add the state $*$ and extend
the partial maps $(-)\cdot a: S\rightharpoonup S$ to the (total) maps
$(-)\cdot a: S\sqcup \{*\} \to S\sqcup \{*\}$ acting by $s\cdot a= *$
if $s\cdot a$ is not defined. Let $S_*= S\sqcup \{*\}$ and $*\cdot a=*$.
The the pair $(M(E,I), S_*)$ consists of a trace monoid with the total action on the set $S_*$.
This pair is called the state space with an augmentation.

For example, the previous state space gives the augmented state space
$$
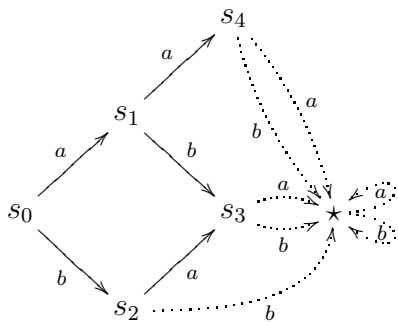
\xymatrix{
& & s_4 \ar@{.>}@/^/[rdd]^a \ar@{.>}@/_/[rdd]_b\\
& s_1 \ar[ru]^{a}\ar[rd]^{b}\\
s_0 \ar[ru]^a \ar[rd]_b && s_3 \ar@{.>}@/^/[r]^a \ar@{.>}@/_/[r]_b & \star
\ar@{.>}@(r,ur)[]^a \ar@{.>}@(r,dr)[]_b\\
& s_2 \ar[ru]_a \ar@{.>}@(r,d)[rru]_b
}
$$

Let $(M(E,I),S)$ be a state space.
Consider an {\em augmented state category} $K_*(S)$ as follows.
Its class of objects is the set  $S_*= S\sqcup \{*\}$.
Morphisms $s\to s'$ are triples $(s,w,s')$
of $s\in S_*$, $s'\in S_*$,
$w\in M(E,I)$.

For any subset $\Sigma\subseteq S_*$, let $K(\Sigma)\subseteq K_*(S)$ denotes
a full subcategory with class of objects $\Sigma$.
For $\Sigma=S$,
$K(S)\subseteq K_*(S)$ will be called a {\em state category}.

\subsection{Homology groups of a small category}

Let $\mC$ be a small category and let $F: \mC\to \Ab$ be a functor
into the category of Abelian groups and homomorphisms.

\begin{definition}\label{defhomolcat}
Let $\mC$ be a small category and let $F: \mC\to \Ab$ be a functor
into the category of Abelian groups and homomorphisms.
Denote by $C_{\diamond}(\mC, F)$ a chain complex of Abelian groups
\begin{displaymath}
 C_n({\mC},F) = \bigoplus_{c_0 \rightarrow \cdots \rightarrow c_n}
F(c_0), \quad n \geq 0,
\end{displaymath}
and homomorphisms
$d_n= \sum\limits_{i=0}^{n}(-1)^i d^n_i:
C_n(\mC,F) \rightarrow C_{n-1}(\mC,F)$, $n>0$,
where
$d^n_i(c_0 \stackrel{\alpha_1}\rightarrow c_1 \stackrel{\alpha_2}\rightarrow
\cdots \stackrel{\alpha_{n}}\rightarrow c_{n}, a)=
$
$$
\left\{
\begin{array}{ll}
(c_1 \stackrel{\alpha_2}\rightarrow \cdots \stackrel{\alpha_n}\rightarrow c_n,
F(c_0\stackrel{\alpha_1}\rightarrow c_1)(a) ) ~~,
& \mbox{if}~ i = 0\\
(c_0 \stackrel{\alpha_1}\rightarrow \cdots \stackrel{\alpha_{i-1}}\rightarrow c_{i-1}
\stackrel{\alpha_{i+1}\alpha_{i}}\rightarrow c_{i+1} \stackrel{\alpha_{i+2}}\rightarrow
\cdots \stackrel{\alpha_n}\rightarrow c_n, a) \quad , & \mbox{if} ~ 1 \leq i \leq n-1\\
(c_0 \stackrel{\alpha_1}\rightarrow \cdots \stackrel{\alpha_{n-1}}\rightarrow c_{n-1},
a ) ~~, & \mbox{if}~ i=n
\end{array}
\right.
$$
For every integer $n\geq 0$, the
{\em $n$-th homology group  $H_n(\mC,F)$ of  $\mC$ with coefficients in $F$}
is the factor groups $\Ker(d_n)/\Imm(d_{n+1})$.
\end{definition}

It is well known that the functors
$H_n(C_{\diamond}(\mC,-)): \Ab^\mC \rightarrow \Ab$
are isomorphic to the left derived functors $\coLim_n^\mC$ of the colimit functor
 $\coLim^\mC: \Ab^\mC \rightarrow \Ab$.

Hence, the Abelian groups $H_n(\mC,F)$ can be defined as homology groups
of the complex
$$
0 \leftarrow \coLim^{\mC}P_0 \leftarrow \coLim^{\mC}P_1 \leftarrow \coLim^{\mC}P_2 \leftarrow \cdots
$$
obtained from a projective resolution
$$
0 \leftarrow F  \leftarrow P_0 \leftarrow P_1 \leftarrow P_2 \leftarrow \cdots
$$
of $F\in \Ab^{\mC}$ by the application of the functor $\coLim^{\mC}$.


\subsection{Homology of state categories, asynchronous systems and Petri nets}


For an arbitrary small category $\mC$,
let $\Delta\ZZ: \mC \to \Ab$ be the functor taking constant values $\ZZ$ at objects and
$1_{\ZZ}: \ZZ \to \ZZ$ at morphisms of $\mC$.

By  \cite{X20042}, an {\em asynchronous system} can be defined as a triple
$(S, s_0, M(E,I))$ where $(S, M(E,I))$ is a state space and
$s_0\in S$ is a distinguished element.
Elements of
$S(s_0)= \{s\cdot\mu| \mu\in M(E,I)\} \subseteq S$ are {\em reachable states}.
{\em Homology groups of asynchronous system with coefficients in
an arbitrary functor  $F: K(S)\to \Ab$}
are Abelian groups
$\coLim_n^{K(S(s_0))}F|_{K(S(s_0))}$.

For a set $B$, denote by $2^B$ the set of all its subsets.

A {\em CE net} \cite{X20042} or {\em Petri net} \cite{win1995}
 is a quintuple $(B,E, pre, post, s_0)$ consisting of finite sets
$B$ and $E$, the maps $pre, post: E\to 2^{B}$,
and a subset $s_0\subseteq B$.

Let ${\mathcal N} =(B,E, pre, post, s_0)$ be a CE net.
Define an relation $I\subseteq E\times E$ as the set of pairs $(a,b)$
for which $(pre(a)\cup post(a))\cap(pre(b)\cup post(b))=\emptyset$.
To every element $e \in E $ we assign a partial mapping $(-)\cdot{e}: 2^B \rightharpoonup 2^B$
defined for $s \subseteq B$ satisfying to the condition
$$
(pre(e)\subseteq s) \quad \& \quad (post(e)\cap s=\emptyset).
$$
In these cases, we take $s\cdot{e}=(s\setminus pre(e)) \cup post(e)$ \cite{maz1987}.
This define a partial action of $M(E,I)$ on the set $2^{B}$.
Assuming $S=2^B$, we get an asynchronous system $(S, s_0, M(E,I))$, which corresponds to
the CE net ${\mathcal N}= (B,E, pre, post, s_0)$.
The homology groups  $H_n({\mathcal N})$ defined as
$\coLim_n^{K(S(s_0))}\Delta\ZZ$ where $S(s_0)$ is the set of reachable states.



In \cite{X20042}, it was built an algorithm for computing the group $H_1(K(S),\Delta\ZZ)$
and hence $H_1({\mathcal N})$. It was formulated the following

\begin{problem}\label{problem1}
Construct an algorithm for computing integral homology of CE nets.
\end{problem}

By the definition of $H_n({\mathcal N})$, this problem will be solved wenn we find
an algorithm to compute the homology groups $H_n(K(S), \Delta\ZZ)$ for the
state categories.

In \cite{X20042}, it was proved that if $M(E,I)$ does not contain triples
of pairwise independent generators, then $H_n(K_*(S),\Delta\ZZ)=0$ for $n>2$.
It was formulated

\begin{problem}\label{problem2}
Let $n>0$ be the maximal number of pairwise independent generators.
Prove that $H_k(K_*(S),F)=0$ for any $k>n$ and for any functor $F: K_*(S)\to \Ab$.
\end{problem}

In the case of finite $E$, this conjecture proved by L. Yu. Polyakova \cite{pol2007}.
In general solved by the author \cite{X20083}.

Problem \ref{problem2} could not be solved for a long time.
We present a way to solve this problem.
Detailed proofs will be published shortly.

\section{Semicubical sets and generalized tori}

Recall the definition of semicubical set and its geometric realization.
Get acquainted with generalized tori and
assign to any partial trace monoid action a  semicubical set.

\subsection{Semicubical sets}

Let $\Box_+$ be the category of posets $\II^n$, $n\in \NN$, where $\II$
is the set $\{0,1\}$ ordered by $0<1$. Morphisms in $\Box_+$ are
increasing maps admitting a decomposition
in the composition of maps $\delta_i^{k,\varepsilon}: \II^{k-1}\rightarrow \II^k$,
$1 \leq i \leq k$, $\varepsilon\in \II$
defined as $\delta_i^{k,\varepsilon}(x_1, \cdots, x_{k-1})=
(x_1, \cdots, x_{i-1}, \varepsilon, x_i, \cdots, x_{k-1}).$

A {\em semicubical set}
is any functor
$X: \Box_{+}^{op} {\rightarrow}\Set$. In \cite{gou1995}, it is called
{\em precubical set}.
Morphisms between semicubical sets are defined as natural transformations.
Any semicubical set can be given by a pair
$(X_n, \partial_i^{n,\varepsilon})$ consisting of
sequence of sets
$(X_n)_{n\in \NN}$ and a family of maps
$\partial_i^{n,\varepsilon}: X_n \rightarrow X_{n-1}$,
defined for $1\leq i\leq n$, $\varepsilon\in \{0,1\}$, and satisfying to the
condition
$$
\partial_i^{n-1,\alpha}\circ \partial_j^{n,\beta} =
\partial_{j-1}^{n-1,\beta}\circ \partial_i^{n,\alpha}~,
\mbox{ for } \alpha,\beta \in \{0,1\}, n\geq 2 \mbox{ and } 1\leq i< j\leq n.
$$
These maps will be equal $\partial_i^{k,\varepsilon}=X(\delta_i^{k,\varepsilon})$.

Semicubical objects in an arbitrary category $\mA$ are defined similarly
as functors $\Box_+^{op}\to \mA$.

\subsection{Geometric realization}

Let $X\in \Set^{\Box_+^{op}}$ be a semicubical set.
Its the {\em geometric realization} \cite{fah2005} is defined as the topological quotient space

$$
|X|_{\Box_+} = \coprod\limits_{n\in \NN} X_n\times [0,1]^n /\equiv
$$
with respect the smallest equivalence relation satisfying
$$
(\partial_i^{n,\nu}x, t_1, \cdots , t_{n-1})\equiv (x, t_1, \cdots, t_{i-1}, \nu, t_i, \cdots , t_{n-1}),
$$
for  all $n\geq 0$, $\in\{0,1\}$, $1\leq i\leq n$, $t_i\in [0,1]$.
Geometric realization determine the functor $|-|_{\Box_+}$ assigning
to every morphism of semicubical sets
 $f: X\to Y$ the continuous map $|f|_{\Box_+}: |X|_{\Box_+}\to |Y|_{\Box_+}$ such that
$|f|_{\Box_+}(x, t_1, \cdots, t_n)=(f(x), t_1, \cdots, t_n)$.
The functor $|-|_{\Box_+}$ can be constructed from the functor $H: \Box_+\to {\rm Top}$,
$H(\II^n)= [0,1]^n$, as in
\cite[Prop. II.1.3]{gab1967}
 by extending to the category of semicubical sets.
 It follows from \cite[Prop. II.1.3]{gab1967} that $|-|_{\Box_+}$
preserves colimits.

\subsection{Generalized tori}

For a trace monoid  $M(E,I)$ with a total order relation $<$ on $E$,
the {\em generalized torus} $T(E,I)$ is
the semicubical set
$(T_n(E,I), \partial^{n,\varepsilon}_i)$ such that
$$
T_n(E,I)=\{(a_1, \cdots, a_n)\in E^n: a_i< a_j ~\& ~(a_i,a_j)\in I \mbox{ for all } 1\leq i<j \leq n\}
$$
and $\partial^{n,\varepsilon}_i(a_1, \cdots, a_n)= (a_1, \cdots, a_{i-1}, a_{i+1}, \cdots, a_n)$,
for all $n\geq 0$, $1\leq i\leq n$, $\varepsilon\in \{0,1\}$.

For example, if $E=\{a_1, \cdots, a_n\}$ ordered by $a_1<a_2<\cdots a_n$
with $I$ consisting of all pairs $(a_i,a_j)$ for which $i\not=j$,
then the geometric realization $|T(E,I)|_{\Box_+}$ is homeomorphic to the usual
$n$-dimensional torus.

\subsection{Semicubical set of a state set}

Let $(M(E,I),S)$ be a state space with a total relation $<$ on $E$.
Assign to it the semicubical set
$Q(E,I,S)$ with
\begin{multline*}
Q_n(E,I,S)=
\{(x, a_1, \cdots, a_n)\in S_*\times T_n(E,I)| \\
a_i <a_j ~\& ~(a_i,a_i)
\mbox{ for all } 1\leq i< j\leq n \}.
\end{multline*}
with the boundary maps
$\partial_i^{n,\varepsilon}(x,a_1, \cdots, a_n)=
(x\cdot a_i^{\varepsilon}, a_1, \cdots, a_{i-1}, a_{i+1}, \cdots, a_n)$
for $1\leq i\leq n$, $n\geq 1$, $\varepsilon\in \{0,1\}$.
Here $a^0=1$ and $a^1=a$.

\begin{example}
Consider the state space consisting of  $S=\{s_0, s_1, s_2, s_3, s_4, s_5\}$,
$E=\{a,b\}$, $I=\{(a,b), (b,a)\}$. Elements in $Tran$ are triples $(s, e, s')$
corresponding to arrows
 $s\stackrel{e}\to s'$ in the following diagram:
$$
\xymatrix{
s_3 \ar[r]^a & s_4 \ar[r]^a & s_5\\
s_0 \ar[r]^a \ar[u]_b & s_1 \ar[r]^a \ar[u]_b & s_2 \ar[u]_b
}
$$

The topological space $|Q(E,I,S)|_{\Box_+}$ can be obtained
from the union of unit squares
$$
\xymatrix{
\star \ar@{-}[r]^a & \star \ar[r]^a & \star\\
\star \ar@{-}[r]^a \ar@{-}[u]_b & \star \ar@{-}[r]^a \ar@{-}[u]_b & \star \ar@{-}[u]_b\\
s_3 \ar@{-}[r]^a \ar@{-}[u]_b & s_4 \ar@{-}[r]^a \ar@{-}[u]_b & s_5 \ar@{-}[r]^a  \ar@{-}[u]_b
							& \star \ar@{-}[r]^a & \star\\
s_0 \ar@{-}[r]^a \ar@{-}[u]_b & s_1 \ar@{-}[r]^a \ar@{-}[u]_b & s_2 \ar@{-}[r]^a  \ar@{-}[u]_b
							& \star \ar@{-}[r]^a \ar@{-}[u]_b & \star\ar@{-}[u]_b
}
$$
by the identifying the vertexes $\star$ with each other,
and by identifying the segments $\xymatrix{\star \ar@{-}[r]^a & \star}$ with each other,
and with similar identifications for the segments
$\xymatrix{\star \ar@{-}[r]^b & \star}$ and squares
$$\xymatrix{\star \ar@{-}[d]_b \ar@{-}[r]^a & \star\ar@{-}[d]^b\\
	\star \ar@{-}[r]^a & \star
}
$$

Geometric realization can be interpreted as the topological space of intermediate states
of computational processes.

\end{example}

\subsection{Homology groups of semicubical sets}

To solve Problems \ref{problem1} and \ref{problem2}, we need an information from the article
\cite{X20081}.

Given a semicubical set  $X\in \Set^{\Box_+^{op}}$,
let
$\Box_+/X$ be the category
with objects
$\sigma\in \coprod\limits_{n\in\NN}X_n$. Its morphisms between
$\sigma\in X_m$ and $\tau\in X_n$ are triples
$(\alpha, \sigma, \tau)$,  $\alpha\in \Box_+(\II^m,\II^n)$,
satisfying
the relation $X(\alpha)(\tau)=\sigma$.
{\em Homological system on a semicubical set
$X$} is an arbitrary functor $F:(\Box_+/X)^{op}\rightarrow \Ab$.

Given a semicubical set $X$ and a homological system $F$, consider
Abelian groups
$C_n(X,F)=\bigoplus\limits_{\sigma\in X_n}F(\sigma)$.
Let $d_i^{n,\varepsilon}: C_n(X,F)\rightarrow C_{n-1}(X,F)$ be the homomorphisms
$$
\bigoplus\limits_{\sigma\in X_n}F(\sigma) \stackrel{d_i^{n,\varepsilon}}\longrightarrow
		\bigoplus\limits_{\sigma\in X_{n-1}}F(\sigma)
$$
defined on the direct summands for
$1\leq i\leq n$,
$\varepsilon \in \II= \{0, 1\}$, $\sigma\in X_n$, $f\in F(\sigma)$ by the equation
$$
d_i^{n.\varepsilon}(\sigma,f)= (\partial_i^{n,\varepsilon}(\sigma),
F(\delta_i^{n,\varepsilon}, \partial^{n,\varepsilon}_i(\sigma),\sigma)(f))\,.
$$

For $n\geq 0$, the {\em homology groups $H_n(X,F)$ of semicubical set $X$
with coefficients in $F$} are defined as homology of the complex
 $C_{\diamond}(X,F)$ consisting of the groups
  $C_n(X,F)= \bigoplus\limits_{\sigma\in X_n}F(\sigma)$
and differentials $d_n=\sum\limits_{i=1}^n (-1)^i (d^{n,1}_i-d^{n,0}_i)$.
Abelian groups $H_n(X,\Delta\ZZ)$ are called to be the $n$th {\em integral}
 homology groups.

\begin{proposition}\label{hsimcube} \cite[Theorem 4.3]{X20081}
For any semicubical set  $X$ and a homological system
$F$ on $X$ there are isomorphisms
$\coLim_n^{(\Box_+/X)^{op}}F  \cong  H_n(X,F)$, for all $n\geq 0$.
\end{proposition}

\begin{proposition}\label{hgeom}
For an arbitrary semicubical set $X$ and integer $n\geq 0$, the group
$H_n(X,\Delta\ZZ)$ is isomorphic to the $n$th singular homology group
of the topological space $|X|_{\Box_+}$.
\end{proposition}

\section{Homology of factorization category}

In this section, we study and apply the Leech homology and cohomology
 groups of trace monoids.

\subsection{Factorization category}

Let $\mC$ be a small category.
Given $\alpha\in \Mor\mC$,
denote by $\cod\alpha$ its codomain and $\dom\alpha$ the domain.

 The {\em factorization category} $Fact(\mC)$
has objects $\Ob(Fact(\mC))=\Mor\mC$,
and for every $\alpha, \beta \in \Mor(\mC)$
each element of $Fact(\mC)(\alpha,\beta)$ is determined by a pair $(f,g)$
of $f,g\in \Mor(\mC)$ making commutative the diagram
$$
\xymatrix{
	\cod\alpha \ar[r]^g & \cod\beta\\
	\dom\alpha\ar[u]^{\alpha}  & \dom\beta\ar[u]_{\beta}\ar[l]^f
}
$$

For example, any monoid $M$ considered as a small category with
unique object has the factorization category $Fact(M)$
such that $\Ob(Fact(M))=M$. Morphisms are given
by quadruples $\alpha\stackrel{(f,g)}\to\beta$ of $f, \alpha,\beta, g \in M$
satisfying $g\alpha f=\beta$.

\subsection{Leech homology of generalized tori}

In this subsection, we present the results published in the articles
 \cite{X20084} and \cite{X2011}.

{\em Leech homology groups of a monoid $M$ with coefficients in a functor}
$F: Fact(M)^{op} \to \Ab$ are defined as the groups $H_n(Fact(M)^{op},F)$, $n\geq 0$.

Given trace monoid $M(E,I)$, let $\mS: \Box_+/T(E,I) \to Fact(M(E,I))$
be the functor assigning to each $(a_1, \cdots, a_n)\in \Ob \Box_+/T(E,I)$
the object $a_1 \cdots a_n \in M(E,I)=\Ob Fact(M(E,I))$.
Each morphism of the category $\Box_+/T(E,I)$ can be
decomposed into a composition of morphisms of the form
$(\delta^{n,\varepsilon}_i), (a_1, \cdots, a_{i-1}, a_{i+1}, \cdots, a_n), (a_1, \cdots, a_n)$.
Therefore, it suffices to define $\mS$ on the morphisms of this kind.
Let
\begin{multline*}
\mS(\delta^{n,\varepsilon}_i, (a_1, \cdots, a_{i-1}, a_{i+1}, \cdots, a_n), (a_1, \cdots, a_n))=\\
(a_1 \cdots a_{i-1} a_{i+1} \cdots a_n \stackrel{(a^{1-\varepsilon}, a^{\varepsilon})}
\longrightarrow
a_1\cdots a_n)
\end{multline*}
where $a^0=1$, and $a^1=a$.

\begin{theorem}\label{faccub}
If $E$ does not contain infinite subsets of pairwise independent elements,
then there are natural in $F\in \Ab^{Fact(M(E,I))^{op}}$ isomorphisms
$$
H_n(Fact(M(E,I))^{op},F) \cong H_n(T(E,I),F\circ\mS^{op}).
$$
\end{theorem}

In the case of a finite set E, this theorem allows us to construct a finite
complex for computing the Leech homology groups.

\subsection{Global dimension of a trace monoid}

Cohomologies of a small categories we define by right derives of the functor $\Lim_{\mC}: \Ab^{\mC}\to \Ab$:

Let $\mC$ be a small category and let $F: \mC\to \Ab$ be a functor.
The category $\Ab^{\mC}$ has enough injectives. Hence
there is an injective resolution
$0\to F \to F^0 \to F^1 \to F^2 \to \cdots$.
Applying the functor $\Lim_{\mC}: \Ab^{\mC}\to \Ab$ to this resolution leads to a complex
$$
	0 \stackrel{d^{-1}}\to \Lim_{\mC}F^0  \stackrel{d^0}\to \Lim_{\mC}F^1
\stackrel{d^1}\to \Lim_{\mC}F^2 \to \cdots
$$
The {\em $n$th cohomology group of $\mC$ with coefficients in $F$} is defined as $H^n(\mC,F)= \Ker d^n/\Imm d^{n-1}$.

Given semicubical set $X$ and a functor $G: \Box_+/X \to \Ab$, define {\em cohomology groups $H^n(X,G)$}
of $X$ with coefficients in $G$
similarly to homology groups of semicubical set. Easy to see that $H^n(X,G)\cong H^n(\Box_+/X,G)$.

The proof of \cite[Theorem 2.2]{X20084} contains the assertion that for each $\alpha\in \Ob Fact(M(E,I))$,
$H_n(\mS/\alpha, \Delta\ZZ)=0$ for $n>0$, and $H_0(\mS/\alpha,\Delta\ZZ)=\ZZ$.
Hence, it follows from the Oberst Theorem \cite[Prop. 1]{X2011}  the following assertion.

\begin{theorem}\label{cohfact}
For any functors $F: Fact(M(E,I))\to Ab$ and for all $n\geq 0$,
there are isomorphisms $H^n(Fact(M(E,I)),F) \cong H^n(T(E,I), F\circ\mS)$.
\end{theorem}

Given Abelian category $\mA$ its {\em global dimension} $\gld\mA$ is a supremum
of $n\geq 0$ for which the functors $\Ext^n(-,=)$ are not equal to $0$.
Theorem \ref{cohfact} leads us to the following generalization
of Hilbert's Syzygy Theorem.

\begin{theorem}\label{globdim}
Let $\mA$ be an Abelian category with coproducts and let $M(E,I)$ be a trace monoid.
If a maximal cardinality of pairwise independent elements of $E$ equals $n<\infty$,
then
$$
	\gld \mA^{M(E,I)} = n + \gld \mA
$$
in each of the following cases:
\begin{enumerate}
\item $\mA$ has exact coproducts (i.e. $\mA$ satisfies to the axiom AB4),
\item $\mA$ has enough projectives.
\end{enumerate}
\end{theorem}
\begin{conjecture}
This is true for all Abelian categories with coproducts.
\end{conjecture}

\begin{example}
Let $k$ be a field and $E= \{x_1, x_2, x_3, x_4, x_5\}$ be the set of variables.
Suppose that the independence relation $I\subset E\times E$
is given by the following graph
with vertexes $E$ and edges $I$:
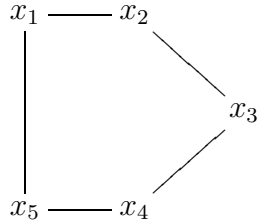
$$
\xymatrix{
x_1 \ar@{-}[r]\ar@{-}[dd] & x_2 \ar@{-}[rd] \\
&& x_3 \\
x_5 \ar@{-}[r] & x_4 \ar@{-}[ru]
}
$$
Denote by $k\langle x_1, x_2, x_3, x_4, x_5\rangle$ the noncommutative polynomial ring in
five variables. Let $(I)$ be the ideal of  $k\langle x_1, x_2, x_3, x_4, x_5\rangle$ generated
by polynomials $x_u x_v - x_v x_u$ for which $(x_u,x_v)\in I$, $1\leq u, v\leq 5$.
The maximal number of independent variables equals $2$. By Theorem \ref{globdim}, we have
$$
\gld k\langle x_1, x_2, x_3, x_4, x_5\rangle/ (I)= 2.
$$
\end{example}

\subsection{Homology of augmented state category}

Consider the functor $\cod: Fact(\mC)\to \mC$, $\alpha\mapsto \cod(\alpha)$,
$(\alpha \stackrel{(f,g)}\longrightarrow \beta) \mapsto$. For any $c\in \Ob\mC$, $H_n(\cod/c,\Delta\ZZ)=0$
for $n$.
\begin{proposition}\label{facus}
Given a small category $\mC$ and a functor $F: \mC^{op}\to \Ab$, there exists
an isomorphisms $\coLim^{\mC^{op}}_n F \cong \coLim^{Fact(\mC)^{op}}_n F\circ\cod^{op}$ for all $n\geq 0$.
\end{proposition}

Given a state space $(M(E,I), S_*)$ and a functor $F: K_*(S)\to \Ab$
there are isomorphisms $H_n(K_*(S),F)\cong H_n(M(E,I)^{op},\overline{F})$
where $\overline{F}=\bigoplus\limits_{x\in S_*}F(x)$ is Abelian group with the right action
$(x,f)\cdot\mu= (x\mu, F(x\stackrel{\mu}\to x\mu)(f))$.
By Proposition \ref{facus} and Theorem \ref{faccub} we obtain the following complex
for the computing the homology of the state space.

\begin{theorem}
If $M(E,I)$ contains no infinite subsets of pairwise independent generators,
then $H_n(K_*(S),F)$ are isomorphic to $n$th homology groups of the complex
\begin{multline*}
0 \leftarrow
\bigoplus\limits_{x\in S_*} F(x) \stackrel{d_1}\leftarrow
\bigoplus\limits_{(x,a_1)\in Q_1(E,I,S)} F(x)
\stackrel{d_2}\leftarrow \bigoplus\limits_{{(x, a_1, a_2)\in Q_2(E,I,S)}
} F(x)
\leftarrow  \cdots \\
\cdots \leftarrow
\bigoplus\limits_{(x, a_1, \cdots, a_{n-1})\in Q_{n-1}(E,I,S)
 }
F(x)
\stackrel{d_n}\longleftarrow
\bigoplus\limits_{{(x, a_1, \cdots, a_n)\in Q_n(E,I,S)} }
F(x) \leftarrow \cdots~,
\end{multline*}
with differentials
\begin{multline*}\label{diff1}
d_n(x,a_1, \cdots, a_n,f) = \\
\sum_{s=1}^n(-1)^s (
(x\cdot a_s, a_1, \cdots, \widehat{a_s}, \cdots, a_n,
F(x\stackrel{a_s}\rightarrow x\cdot a_s)(f))\\
- (x, a_1, \cdots, \widehat{a_s}, \cdots, a_n, f) )
\end{multline*}
\end{theorem}

So, we have the following solution of Problem \ref{problem2}.

\begin{corollary}
If the cardinality of pairwise generators of $M(E,I)$ not greater than $n$,
then $H_k(K_*(S),F)=0$ for all $k>n$.
\end{corollary}

In addition, we have a complex of finitely generated abelian groups for calculating
the integral homology $H_n(K_*(S),\Delta\ZZ)$ of augmented state category.

\begin{example}
Consider a state space
$\Sigma=(S, E, I, Trans)$, $S=\{s_0, s_1\}$, $E=\{a, b\}$, $I=\{(a,b), (b,a)\}$,
$Tran= \{(s_0, a, s_0), (s_0, b, s_1), (s_1,a,s_1)\}$. The set consists of two elements
with the partial action of the free commutative monoid generated by $a$ and $b$.
Let us calculate the groups  $H_n(K_*(S), \Delta\ZZ)$.

We add the state $\star$
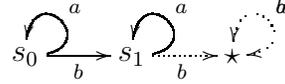
$$
\xymatrix{
	s_0 \ar@(r, u)_{a} \ar[r]_b & s_1 \ar@(r,u)_a \ar@{.>}[r]_b & \star\ar@{.>}@(r,u)_a \ar@{.>}@(u, r)^b
}
$$
and write down the matrixes of differentials.
Since $|S_*|= 3$, $|Q_1(E,I,S_*)|=6$, $|Q_2(E,I,S_*)|=3$, the complex consists of Abelian groups
$$
0 \leftarrow \ZZ^3 \stackrel{d_1}\leftarrow \ZZ^6 \stackrel{d_2}\leftarrow \ZZ^3 \leftarrow 0
$$
The differential $d_1(s,e)= -s\cdot e + s$ is defined by the matrix:

\begin{gather*}
\quad
\begin{array}{cccccccc}
 & (s_0,a) & ~~(s_0,b) & (s_1,a) & ~~(s_1,b) & ~~~(*,a) & ~~~~(*,b)
\end{array}
\\
\begin{array}{l}
s_0\\
s_1\\
\star
\end{array}
\left(
\begin{array}{cccccc}
+1-1 & ~~+1 & ~~~0 & ~~~~0& ~~~~0& ~~~~0\\
0    & ~~-1 & ~~~+1-1 & ~~~~+1& ~~~~0& ~~~~0\\
0    & ~~~0 & ~~~0 & ~~~~-1& ~~-1+1& ~~-1+1
\end{array}
\right)
\end{gather*}
The differential $d_2(s, e_1, e_2)=-(s*e_1, e_2)+(s,e_2)+(s*e_2,e_1)-(s,e_1)$ has
the matrix:
\begin{gather*}
\quad
\begin{array}{cccc}
 & ~~~~~ (s_0,a, b) & (s_1,a, b) & (\star, a, b)
\end{array}
\\
\begin{array}{l}
(s_0,a)\\
(s_0,b)\\
(s_1,a)\\
(s_1,b)\\
(\star,a)\\
(\star,b)
\end{array}
\left(
\begin{array}{ccc}
-1 & 0 & 0\\
-1+1 & 0 & 0\\
+1 & -1 & 0\\
0 & -1+1 & 0\\
0 & +1 & +1-1\\
0 & 0 & -1+1
\end{array}
\right)
\end{gather*}
Using reduction of these matrices to Smith normal form, we obtain
$H_0(K_*(S), \Delta\ZZ)=\ZZ$, $H_1(K_*(S), \Delta\ZZ)=\ZZ^2$,
$H_2(K_*(S), \Delta\ZZ)=\ZZ^1$, and $H_n(K_*(S), \Delta\ZZ)=0$
for all $n\geq 3$.
\end{example}

\subsection{Homology of Mazurkiewicz trace languages}

Given $v,w \in M(E,I)$, we let $v\leq w$ if there exists $u\in M(E,I)$ such that $vu=w$.
This relation makes $M(E,I)$ into a partially ordered set,
which we denote by $P(E,I)$.
A {\em trace language} is any set of traces.
\begin{definition}
A set $L\subseteq M(E,I)$ is {\em prefix closed} if
for all $v\in M(E,I)$ and $w\in L$ the relation $v<w$ implies $v\in L$.
\end{definition}

Let $L\subseteq M(E,I)$ be a prefix closed trace language.
We have the pair $(M(E,I),L)$ consisting of the trace monoid
with the following partial action for $v\in L$, $\mu\in M(E,I)$.
$$
v\cdot\mu=\left\{
\begin{array}{ll}
v\mu, & \mbox{if} ~~ v\mu\in L\\
undefined, & otherwise.
\end{array}
\right.
$$

For any functor $F: K_*(L)\to \Ab$, we can consider the homology groups
$H_n(K_*(L), F)$.
The groups $H_n(K_*(L), \Delta\ZZ)$ are called
{\em integral homology groups}.

\subsection{Homology groups of the poset of traces}

Given prefix closed language $L\subseteq M(E,I)$,
let $\ZZ[L]: P(E,I)\to \Ab$ be a functor with values
$\ZZ[L](v)=\ZZ$ for $v\in L$ and $\ZZ[L](v)=0$, otherwise.
For $u\leq v\in L$, we will define  $\ZZ[L](u\leq v)=1_{\ZZ}$.
We study the homology groups $H_n(P(E,I), \ZZ[L])$ of the poset
$P(E,I)$ and their relationship with  $H_n(K_*(L), \Delta\ZZ)$.

Denote by $p_n$ the cardinality of the set of $n$-cliques in the
graph $(E,I)$. In particular, $p_0=1$ as the number of empty subsets in $E$,
$p_1=|E|$.
For example, if $(E,I)$ is the graph
$$
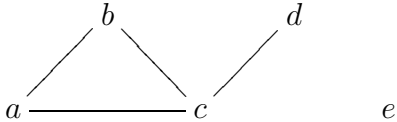
\xymatrix{
& b && d\\
a\ar@{-}[ru] \ar@{-}[rr]&& c\ar@{-}[lu] \ar@{-}[ru] && e
}
$$
then $p_0=1$, $p_1=5$, $p_2=4$, $p_3=1$.

\begin{theorem}
$H_n(K_*(L), \Delta\ZZ)\cong H_n(P(E,I),\ZZ[L])\oplus \ZZ^{(p_n)}$.
\end{theorem}

Given a partially ordered set $P$, let $\widetilde{H_n}(P)$
be the reduced singular homology of the classifying space $B(P)$.
It is not hard to see that $H_n(P(E,I),\ZZ[L])\cong \widetilde{H}_{n-1}(P(E,I)\setminus L)$
for $n\geq 1$.
\begin{corollary}
$H_n(K_*(L), \Delta\ZZ)\cong \widetilde{H}_{n-1}(P(E,I)\setminus L)\oplus \ZZ^{(p_n)}$
for all $n\geq 1$.
\end{corollary}
We see that $H_1(K_*(L),\Delta\ZZ)$ is a free Abelian group.
\begin{conjecture}
For any trace monoid $M(E,I)$ with partial action on a set $S$, the Abelian group
$H_1(K_*(S), \Delta\ZZ)$
is free.
\end{conjecture}

Note the following homological properties of partially ordered set of traces.
We assume that the language of traces $L$ is prefix closed.
\begin{itemize}
\item If $I=\{(a,b)\in E\times E| a\not=b\}$ and hence $M(E,I)$ is commutative, then $H_n(P(E,I),\ZZ[L])=0$ for all $n\geq 1$.
\item If $I=\emptyset$ and hence $M(E,I)$ is free,
then $H_n(P(E,I),\ZZ[L])=0$ for all $n\geq 2$.
\item For arbitrary finitely generated Abelian groups $A_1$, $A_2$, ..., $A_n$
with free $A_1$, there exists a trace monoid $M(E,I)$ such that
$H_n(P(E,I), \ZZ[\{1\}])\cong A_k$ for all $1\leq k\leq n$.
\end{itemize}


\subsection{Baues-Wirsching homology of the state category}

Let $M(E,I)$ be an arbitrary trace monoid an let $X$ be a right $M(E,I)$-set.
Recall that $K(X)$ denotes the category of states with objects $x\in X$ and
morphisms $x\stackrel{\mu}\to x\mu$ for $x\in X$ and $\mu\in M(E,I)$.
Considering $M(E,I)$ as a category with an unique object
we can define a functor $U: K(X)\to M(E,I)$ assigning to each morphism
$x\stackrel{\mu}\to x\mu$ the morphism $\mu\in M(E,I)$.
Applying the functor $Fact$ to $U$, we can consider a functor
$Fact(U): Fact(K(X))\to Fact(M(E,I))$.
For any functor $F: Fact(K(X))^{op}\to \Ab$, there exists its Kan extension
$Lan^{Fact(U)^{op}}: Fact(K(M(E,I)))^{op}\to \Ab$ \cite{mac1998}.
\begin{theorem}
Given functor $F: Fact(K(X))^{op}\to \Ab$, there exist isomorphisms
$$
  H_n(Fact(K(X))^{op}, F)\cong H_n(Fact(M(E,I))^{op}, Lan^{Fact(U)^{op}}F)
$$
for all $n\geq 0$.
\end{theorem}

\subsection{The solution of Problem 1}

Recall that a state space $(M(E,I),S)$ is a trace monoid with a partial
action on $S$. The category of states $K(S)\subset K_*(S)$ is the full subcategory
with objects $s\in S$. Denote by $\ZZ S$ the free Abelian group generated by $s\in S$.
Let $\overline{Q}_n(E,I,S)=\{(s, a_1, \cdots, a_n)\in S\times T_n(E,I)|
s a_1\cdots a_n\not=\star\}$.
\begin{theorem}
Given a state space $(M(E,I),S)$, the groups $H_n(K(S),\Delta\ZZ)$
are isomorphic to the homology groups of the complex
\begin{multline*}
0 \leftarrow \ZZ(S) \stackrel{d_1}\leftarrow \ZZ \overline{Q}_1(S,E,I)
\stackrel{d_2}\leftarrow
\ZZ \overline{Q}_2(S,E,I)
\leftarrow \cdots \\
\cdots \leftarrow \ZZ \overline{Q}_{n-1}(S,E,I) \stackrel{d_n}\leftarrow
\ZZ \overline{Q}_n(S,E,I)
\leftarrow \cdots
\end{multline*}
with differentials
\begin{multline*}
d_n(s, a_1, \cdots, a_n)= \sum\limits_{i=1}^n(-1)^i(sa_i, a_1, \cdots, a_{i-1}, a_{i+1}, \cdots, a_n)\\
			- \sum\limits_{i=1}^n(-1)^i(s, a_1, \cdots, a_{i-1}, a_{i+1}, \cdots, a_n)
\end{multline*}
\end{theorem}

Consider an example of computing the homomology groups of a state category.
\begin{example}
Let $M(E,I)$ be a commutative trace monoid
generated by two elements and let $S$ consists of two elements.
That is $E=\{a,b\}$, $I=\{(a,b), (b,a)\}$, $S=\{s_0, s_1\}$. The generators act by
$s_0 a= s_0$, $s_0 b=s_1$, $s_1 a= s_1$ as it is shown in the following picture.
$$
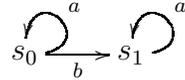
\xymatrix{
	s_0 \ar@(r, u)_{a} \ar[r]_b & s_1 \ar@(r,u)_a
}
$$
The complex consists of abelian groups
$$
C_0=\ZZ\{s_0, s_1\},\quad C_1=\ZZ\{(s_0, a), (s_0,b), (s_1, a)\}, \quad
C_2=\ZZ\{(s_0, a, b)\}.
$$
We have a complex  $0 \leftarrow \ZZ^2 \stackrel{d_1}\leftarrow \ZZ^3 \stackrel{d_2}\leftarrow \ZZ \leftarrow
0 \leftarrow 0 \leftarrow \cdots$.
The differential $d_1$ is described by the following matrix.

\begin{gather*}
\quad
\begin{array}{cccc}
 & (s_0,a) & ~~(s_0,b) & (s_1,a)
\end{array}
\\
\begin{array}{l}
s_0\\
s_1
\end{array}
\left(
\begin{array}{cccccc}
1-1 & ~~ 1 & ~~~0 \\
0    & ~~-1 & ~~~1-1
\end{array}
\right)
\end{gather*}

The differential $d_2$ has
the following matrix.
\begin{gather*}
\quad
\begin{array}{cccc}
 & ~~~~~ (s_0,a, b)
\end{array}
\\
\begin{array}{l}
(s_0,a)\\
(s_0,b)\\
(s_1,a)
\end{array}
\left(
\begin{array}{c}
-1 \\
-1+1\\
+1
\end{array}
\right)
\end{gather*}
Using the reduction to the Smith normal forms, we get
$$
H_0(K(S),\Delta\ZZ)=\ZZ, ~ H_1(K(S), \Delta\ZZ)=\ZZ, ~ H_n(K(S),\Delta\ZZ)=0 ~ \mbox{for all}~ n\geq 2.
$$

\end{example}

\subsection{Homology groups of CE nets}

For the computing the homology groups of a finite CE net,
we first construct the state space $(M(E,I), S(s_0))$.
Then we can compute $H_n(K(S(s_0)),\Delta\ZZ)$ by the method described above.

Let, for example, $\mN$ be the following CE net.
$$
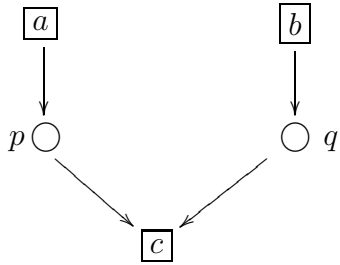
\xymatrix{
~~\boxed{a}\ar@<1ex>[d] & & \boxed{b}\ar[d] \\
~p\,{\bigcirc} \ar[rd] && ~~\bigcirc\ar[ld] {\,q}\\
  & \boxed{c}\\
}
$$
The corresponding trace monoid $M(E,I)$ is defined by
$E=\{a,b,c\}$ and $I=\{(a,b),(b,a)\}$.
The set of states $S$ consists of all subsets
$s\subseteq \{p,q\}$.
The corresponding asynchronous system  $(M(E,I),S, s_0)$ is defined
by $s_0=\emptyset$ and a partial action of $M(E,I)$ shown in the following figure.
$$
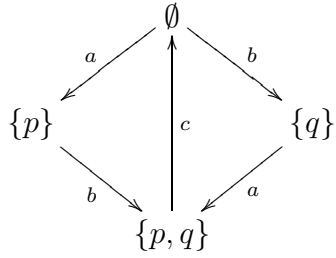
\xymatrix{
& \emptyset \ar[ld]_a \ar[rd]^b\\
\{p\}\ar[rd]_b && \{q\}\ar[ld]^a\\
& \{p,q\}\ar[uu]_c
}
$$
That is $\emptyset\cdot a=\{p\}$,  $\emptyset\cdot b=\{q\}$,  $\{p\}\cdot b=\{p,q\}$,
$\{q\}\cdot a=\{p,q\}$, and $\{p,q\}\cdot c = \emptyset$.
All states are admissible. Hence $S(s_0)=S$.
The complex consists of the Abelian groups
\begin{gather*}
C_0=\ZZ\{\emptyset, \{p\}, \{q\}, \{p,q\}\}\cong\ZZ^4,\\
C_1=\ZZ\{(\emptyset, a), (\emptyset, b),
(\{p\},b), (\{q\},a), (\{p, q\},c) \}\cong \ZZ^5,\\
C_2= \ZZ\{(\emptyset,a,b)\}\cong \ZZ.
\end{gather*}

The differential $d_1(s,e)=-s\cdot e+ s$ has the following matrix.

\begin{gather*}
\quad
\begin{array}{ccccccc}
\quad &~~~ \quad  (\emptyset,a) & \quad (\emptyset,b) & (\{p\},b) & (\{q\},a) & (\{p,q\},c)
\end{array}
\\
\begin{array}{c}
\emptyset\\
\{p\}\\
\{q\}\\
\{p,q\}
\end{array}
\left(
\begin{array}{cccccc}
~~ ~1    & \quad~~~1 &\quad ~~~0 & \quad~~~~0& \quad~~~~0&\\
~~ -1    &\quad ~~~0 &\quad ~~~1 &\quad ~~~~1& \quad~~~-1\\
~~ ~0    &\quad ~~-1 &\quad ~~~0 &\quad ~~~~0& \quad~~~~0\\
~~ ~0    &\quad ~~~0 &\quad ~~-1 &\quad ~~~~-1& \quad~~~~1
\end{array}
\right)
\end{gather*}
We have $d_2(\emptyset, a, b)=-(\emptyset\cdot a, b)+(\emptyset, b)+(\emptyset\cdot b, a)
-(\emptyset,a)$. Hence, the matrix of $d_2$ is described
by the matrix
\begin{gather*}
\quad
\begin{array}{cc}
\quad & \qquad (\emptyset,a, b)
\end{array}
\\
\begin{array}{c}
(\emptyset,a)\\
(\emptyset,b)\\
(\{p\},b)\\
(\{q\},a)\\
(\{p,q\}, c)
\end{array}
\left(
\begin{array}{c}
-1 \\
1\\
-1\\
1 \\
0
\end{array}
\right)
\end{gather*}
We have the following complex for the computing $H_n(\mN)$ for all $n\geq 0$.
$$
0 \leftarrow \ZZ^4 \stackrel{d_1}\leftarrow \ZZ^5 \stackrel{d_2}\leftarrow \ZZ
\leftarrow 0 \leftarrow 0\leftarrow \cdots
$$
Using the Smith normal forms, we get
$H_0(\mN)=\ZZ$, $H_1(\mN)=\ZZ$, and
$H_n(\mN)=0$,
for all $n\geq 2$.

\section{Conclusion}

The author believes that the results will help in investigation the Goubault homology 
of asynchronous systems as the homology groups
$H_n(K(S), \ZZ^{\varepsilon})$, $\varepsilon\in \{0,1\}$, 
with coefficients in some suitable systems of Abelian groups. 
You can explore the $n$-deadlocks for asynchronous systems. 
It is possible to find homological signs for the existence of bisimilar equivalence
between asynchronous systems, Petri nets, and trace languages.

\end{document}